\documentclass[11pt]{article}
\usepackage{amsfonts}
\usepackage{graphics}
\usepackage{indentfirst}
\usepackage{cite}
\usepackage{latexsym}
\usepackage[paper=a4paper, left=2.1cm, right=2.1cm, top=2.2cm, bottom=1.5cm, headheight=5.5pt, footskip=0.8cm, footnotesep=0.8cm, centering, includefoot]{geometry}
\usepackage{amsmath}
\allowdisplaybreaks
\usepackage{amssymb}
\usepackage[dvips]{epsfig}
\usepackage{amscd}

\newtheorem{theorem}{Theorem}[section]
\newtheorem{remark}{Remark}[section]

\newtheorem{definition}{Definition}[section]
\newtheorem{lemma}[theorem]{Lemma}

\newtheorem{proposition}[theorem]{Proposition}

\DeclareMathOperator{\loc}{loc}

\DeclareMathOperator{\divv}{div}

\def\norm[#1]#2{\|#2\|_{#1}}

\makeatletter
\@addtoreset{equation}{section}
\makeatother

\title{Strong solutions to the Cauchy problem of the two-dimensional non-baratropic non-resistive magnetohydrodynamic equations with zero heat conduction
\thanks{Supported by China Postdoctoral Science Foundation (No. 2017M610579), the Postdoctoral Science Foundation of Chongqing (No. xm2017015), Fundamental Research Funds for the Central Universities (No. XDJK2017C050), and the Doctoral Fund of Southwest University (No. SWU116033).
}
}
\date{}
\author{Xin Zhong\thanks{School of Mathematics and Statistics,
Southwest University, Chongqing 400715, People's Republic of China ({\tt xzhong1014@amss.ac.cn}).
 }
 }

\begin{document}
\maketitle

\begin{abstract}
This paper concerns the Cauchy problem of the non-baratropic non-resistive  magnetohydrodynamic (MHD) equations with zero heat conduction on the whole two-dimensional (2D) space with vacuum as far field density.
By delicate weighted energy estimates, we prove that there exists a local strong solution provided the initial density and the initial magnetic decay not too slow at infinity.
\end{abstract}


Math Subject Classification: 76W05; 35B65

\section{Introduction}
Let $\Omega\subset\mathbb{R}^n\ (n=2,3)$ be a domain, the motion of a viscous, compressible, and heat conducting magnetohydrodynamic (MHD) flow in $\Omega$ can be described by the full compressible MHD equations
\begin{align}\label{1.1}
\begin{cases}
\rho_{t}+\divv(\rho\mathbf{u})=0,\\
(\rho\mathbf{u})_{t}+\divv(\rho\mathbf{u}\otimes\mathbf{u})
-\mu\Delta\mathbf{u}
-(\lambda+\mu)\nabla\divv\mathbf{u}+\nabla P=
\mathbf{b}\cdot\nabla\mathbf{b}-\frac12\nabla|\mathbf{b}|^2,\\
c_{\nu}[(\rho\theta)_{t}+\divv(\rho\mathbf{u}\theta)]
+P\divv\mathbf{u}-\kappa\Delta\theta
=2\mu|\mathfrak{D}(\mathbf{u})|^2+\lambda(\divv\mathbf{u})^2
+\nu|\nabla\times\mathbf{b}|^2,\\
\mathbf{b}_{t}-\mathbf{b}\cdot\nabla\mathbf{u}
+\mathbf{u}\cdot\nabla\mathbf{b}+\mathbf{b}\divv\mathbf{u}=\nu\Delta\mathbf{b},\\
\divv \mathbf{b}=0.
\end{cases}
\end{align}
Here, $t\geq0$ is the time, $x\in\Omega$ is the spatial coordinate, and $\rho, \mathbf{u}, P=A\rho\theta\ (A>0), \theta, \mathbf{b}$ are the fluid density, velocity, pressure, absolute temperature, and the magnetic field respectively; $\mathfrak{D}(\mathbf{u})$ denotes the deformation tensor given by
\begin{equation*}
\mathfrak{D}(\mathbf{u})=\frac{1}{2}(\nabla\mathbf{u}+(\nabla\mathbf{u})^{tr}).
\end{equation*}
The constant viscosity coefficients $\mu$ and $\lambda$ satisfy the physical restrictions
\begin{equation}\label{1.2}
\mu>0,\ 2\mu+n\lambda\geq0.
\end{equation}
Positive constants $c_\nu,\ \kappa$, and $\nu$ are respectively the heat capacity, the ratio of the heat conductivity coefficient over the heat capacity, and the magnetic diffusive coefficient.

There is huge literature on the studies about the theory of
well-posedness of solutions to the compressible MHD system due to the physical importance, complexity, rich phenomena and mathematical challenges.
However, many physical important and mathematical fundamental problems are still open due to the lack of smoothing mechanism and the strong nonlinearity. The local existence and uniqueness of strong solutions were obtained
by Vol'pert-Hudjaev \cite{VH1972} as the initial density is strictly positive. Kawashima \cite{K1983} first obtained the global existence and
uniqueness of classical solutions when the initial data are close to a non-vacuum equilibrium in $H^3$-norm.
When the initial density allows vacuum, the local well-posedness of strong solutions to the initial boundary value problem of 3D nonisentropic MHD equations has been obtained by Fan-Yu \cite{FY2009}. For general large initial data, Hu-Wang \cite{HW2008,HW2010} proved the global existence of weak solutions with finite energy in Lions' framework for compressible Navier-Stokes equations \cite{L1998,F2004} provided the adiabatic exponent is suitably large.
Recently, Li-Xu-Zhang \cite{LXZ2013} established the global existence and uniqueness of classical solutions to the Cauchy problem for the isentropic compressible MHD system in 3D with smooth initial data which are of small energy but possibly large oscillations and vacuum, which generalized the result for compressible Navier-Stokes equations obtained by Huang-Li-Xin \cite{HLX2012}. Very recently, Hong-Hou-Peng-Zhu \cite{HHPZ2017} improved the result in \cite{LXZ2013} to allow the initial energy large as long as the adiabatic exponent is close to 1 and $\nu$ is suitably large.

For two dimensional problems, when the far-field density is away from vacuum, for any $p\in[2,\infty)$, the $L^p$-norm of a function $\mathbf{u}$ can be bounded by $\|\sqrt{\rho}\mathbf{u}\|_{L^2}$ and $\|\nabla\mathbf{u}\|_{L^2}$. Hence, we can modify the techniques in \cite{L2012} and obtain local existence of strong solutions to the Cauchy problem of MHD equations. However, when
the far-field density equals zero, it seems difficult to bound the $L^p$-norm of $\mathbf{u}$ by $\|\sqrt{\rho}\mathbf{u}\|_{L^2}$ and $\|\nabla\mathbf{u}\|_{L^2}$ for any $p\geq1$, so even the local existence of strong solutions to the Cauchy problem is much more subtle. Recently, by using delicate weighted energy estimates, Li-Liang \cite{LL2014} established the local existence of unique strong solutions for the two-dimensional Cauchy problem of the barotropic compressible Navier-Stokes equations with vacuum as far-field density, then L{\"u}-Huang \cite{LH2015}
extended this result to the 2D compressible MHD equations.
Later on, Lu-Huang \cite{LH2016} proved the local existence of strong solutions to the Cauchy problem of 2D nonisentropic MHD equations with zero heat conduction, which generalized the result of the nonisentropic Navier-Stokes equations obtained by Liang-Shi \cite{LS2015}.

As is well-known, when there is a diffusion term $\Delta\mathbf{b}$ in the magnetic equations, due to the similarity of the momentum equations and the magnetic equations, the study for MHD system has been along with that for the Navier-Stokes one.
By contrast, for the MHD system with no diffusion for the magnetic field,
new difficulties will arise in the study of strong (or classical) solvability to the MHD system because of the lack of smoothing effect on the magnetic equations. Thus it is a natural and interesting problem to investigate solvability of the compressible non-resistive MHD system. Recently, Zhu \cite{Z2015} showed the existence of unique local classical solutions of three-dimensional compressible isentropic magnetohydrodynamic equations with infinite electric conductivity when the initial data are arbitrarily large, contain vacuum and satisfy some suitable conditions. Chen-Zang \cite{CZ2017} proved that the two-dimensional case admits a unique local strong solution provided the initial density and the magnetic field decay not too slow at infinity. Moreover, if the initial data satisfy some additional regularity and compatibility conditions, the strong solution becomes a classical one. Very recently, motivated by \cite{LL2014,CZ2017}, Zhong \cite{Z2017} investigated the local existence of unique strong solutions of 2D Cauchy problem to the compressible non-resistive MHD equations with variable viscosity.
The aim of this paper is to show the local existence of strong solutions to the 2D Cauchy problem \eqref{1.1} with $\kappa=\nu=0$, which will generalize the local theory of
\cite{CZ2017} to the nonisentropic case.

Let $\Omega=\mathbb{R}^2$ and $\kappa=\nu=0$. Without loss of generality, we take $c_\nu=A=1$, then the system \eqref{1.1} can be written as
\begin{align}\label{mhd}
\begin{cases}
\rho_{t}+\divv(\rho\mathbf{u})=0,\\
(\rho\mathbf{u})_{t}+\divv(\rho\mathbf{u}\otimes\mathbf{u})
-\mu\Delta\mathbf{u}
-(\lambda+\mu)\nabla\divv\mathbf{u}+\nabla P=
\mathbf{b}\cdot\nabla\mathbf{b}-\frac12\nabla|\mathbf{b}|^2,\\
P_{t}+\divv(P\mathbf{u})
+P\divv\mathbf{u}=2\mu|\mathfrak{D}(\mathbf{u})|^2+\lambda(\divv\mathbf{u})^2,\\
\mathbf{b}_{t}-\mathbf{b}\cdot\nabla\mathbf{u}
+\mathbf{u}\cdot\nabla\mathbf{b}+\mathbf{b}\divv\mathbf{u}=\mathbf{0},\\
\divv \mathbf{b}=0.
\end{cases}
\end{align}
The present paper is aimed at establishing strong solutions to the Cauchy problem  of the system \eqref{mhd} with the initial condition
\begin{equation}\label{n3}
(\rho,\rho\mathbf{u},P,\mathbf{b})(x,0)=(\rho_0,\rho_0\mathbf{u}_0,P_0,\mathbf{b}_0)(x),\ \ x\in\mathbb{R}^2,
\end{equation}
and the far field behavior
\begin{equation}\label{n4}
(\rho,\mathbf{u},P,\mathbf{b})(x,t)\rightarrow(0, \mathbf{0},0, \mathbf{0}),\ \text{as}\ |x|\rightarrow+\infty,\ t>0.
\end{equation}

Before stating our main result, we first explain the notations and conventions used throughout this paper. For $r>0$, set
\begin{equation*}
B_r \triangleq\left.\left\{x\in\mathbb{R}^2\right|\,|x|<r \right\}.
\end{equation*}
For $1\leq p\leq\infty$ and integer $k\geq0$, the standard Sobolev spaces are denoted by:
\begin{equation*}
L^p=L^p(\mathbb{R}^2),\ W^{k,p}=W^{k,p}(\mathbb{R}^2), \ H^{k}=H^{k,2}(\mathbb{R}^2),
\ D^{k,p}=\{u\in L_{\loc}^1|\nabla^k u\in L^p\}.
\end{equation*}
Now we define precisely what we mean by strong solutions to the problem \eqref{mhd}--\eqref{n4}.
\begin{definition}[Strong solutions]\label{def1}
$(\rho,\mathbf{u},P,\mathbf{b})$ is called a strong solution to \eqref{mhd}--\eqref{n4} in $\mathbb{R}^2\times(0,T)$ if all derivatives involved in \eqref{mhd}--\eqref{n3}
for $(\rho,\mathbf{u},P,\mathbf{b})$ are regular distributions,
and $(\rho,\mathbf{u},P,\mathbf{b})$ satisfies both \eqref{mhd} almost everywhere in $\mathbb{R}^2\times(0,T)$ and \eqref{n3} almost everywhere in $\mathbb{R}^2$.
\end{definition}

Our main result reads as follows:
\begin{theorem}\label{t1}
Let $\eta_0 $ be a positive constant and
\begin{equation}\label{1.6}
\bar x\triangleq\left(e+|x|^2\right)^{\frac12}\log^{1+\eta_0}\left(e+|x|^2\right).
\end{equation}
For constants $q>2$ and $a>1$, assume that the initial data $(\rho_0\geq0, \mathbf{u}_0, P_0\geq0, \mathbf{b}_0)$ satisfy
\begin{align}\label{1.7}
\begin{cases}
\rho_0\bar x^a\in L^1 \cap H^1\cap W^{1,q},\,\,\mathbf{b}_0\bar x^{a}\in L^2\cap L^{q},\ \nabla\mathbf{b}_0\in L^2\cap L^{q},
\\
P_0\in L^1 \cap H^1\cap W^{1,q},\,\,\sqrt{\rho_0}\mathbf{u}_0\in L^2,\,\,\nabla\mathbf{u}_0\in  L^2,\,\,\divv\mathbf{b}_0=0,
\end{cases}
\end{align}
and the following compatibility condition
\begin{align}\label{1.8}
-\mu\Delta\mathbf{u}_0-(\mu+\lambda)\nabla\divv\mathbf{u}_0
+\nabla P_0-\mathbf{b}_0\cdot\nabla\mathbf{b}_0+\frac12\nabla|\mathbf{b}_0|^2
=\sqrt{\rho_0}\mathbf{g}
\end{align}
for some $\mathbf{g}\in L^2(\mathbb{R}^2)$.
Then there exists a positive time $T_0>0$ such that the problem  \eqref{mhd}--\eqref{n4} has a strong solution $(\rho,\mathbf{u},P,\mathbf{b})$ on $\mathbb{R}^2\times (0,T_0]$ satisfying
\begin{align}\label{1.9}
\begin{cases}
\rho\geq0,\ \rho\in C([0,T_0];L^1\cap H^1\cap W^{1,q}),\\
\rho\bar{x}^a\in L^\infty(0,T_0;L^1\cap H^1\cap W^{1,q}),\\
\sqrt{\rho}\mathbf{u},\nabla\mathbf{u},\sqrt{\rho}\dot{\mathbf{u}}\in L^\infty(0,T_0;L^{2}),\
\nabla\mathbf{u}\in L^2(0,T_0;H^1\cap W^{1,q}),\\
\mathbf{b}\in C([0,T_0];L^2\cap L^{q}),\
\mathbf{b}\bar{x}^a\in L^\infty(0,T_0;L^2\cap L^{q}), \\
P\geq0,\ P\in C([0,T_0];L^1\cap H^1\cap W^{1,q}),
\end{cases}
\end{align}
and
\begin{equation}\label{1.10}
\inf\limits_{0\le t\le T_0}\int_{B_{N}}\rho(x,t)dx\ge \frac14\int_{\mathbb{R}^2} \rho_0(x)dx
\end{equation}
for some constant $N >0$.
\end{theorem}

\begin{remark}
When there is no electromagnetic field effect, that is $\mathbf{b}=\mathbf{0}$,
\eqref{mhd} turns to be the full compressible Navier-Stokes equations without heat conduction, and Theorem \ref{t1} is similar to the result of \cite{LS2015}.
Roughly speaking, we generalize the result of \cite{LS2015} to the full compressible non-resistive MHD system with zero heat conduction.
\end{remark}

\begin{remark}
Due to the insufficient smoothness and the absence of the compatibility conditions on the initial data, and the presence of vacuum, for the strong solution $(\rho,\mathbf{u},P,\mathbf{b})$ obtained in Theorem \ref{t1}, the quantity $\mathbf{u}_t$, viewed as a vector valued function on the time interval $(0, T_0)$, is not generally integrable on $(0, T_0)$. As a result, one cannot expect the continuity of $\mathbf{u}$ up to the initial time. However, it follows from \eqref{1.9} and Aubin-Lions lemma, we see that the velocity $\mathbf{u}$ is continuous with respect to $t$ as long as $t>0$. Nevertheless, we can get the continuity of $\rho\mathbf{u}$ at the initial time. This is the main reason that we impose the initial condition on $\rho\mathbf{u}$ instead of $\mathbf{u}$ in \eqref{n3}.

Indeed, by \eqref{1.9}, we immediately have
\begin{align}\label{10}
& \rho\mathbf{u}=\sqrt{\rho}\cdot\sqrt{\rho}\mathbf{u}\in L^\infty_tL^2_x,\\
& \rho\nabla\mathbf{u}\in L^\infty_tL^2_x\label{20},\\
& \nabla\rho\bar{x}^a\in L^\infty_tL^q_x\ \ \text{for}\ q>2\label{30}.
\end{align}
By \eqref{2.6} in the next section and \eqref{1.9}, we derive that
\begin{equation*}
\mathbf{u}\bar{x}^{-a}\in L^\infty_tL^p_x\ \ \text{for any}\ p>2,
\end{equation*}
which combined with \eqref{30} and H{\"o}lder's inequality leads to
\begin{equation}\label{40}
\nabla\rho\cdot\mathbf{u}=\nabla\rho\bar{x}^a\cdot
\mathbf{u}\bar{x}^{-a}\in L^\infty_tL^2_x.
\end{equation}
Thus, we infer from \eqref{20} and \eqref{40} that
\begin{equation*}
\nabla(\rho\mathbf{u})=\rho\nabla\mathbf{u}+\nabla\rho\cdot\mathbf{u}\in L^\infty_tL^2_x,
\end{equation*}
which along with \eqref{10} yields
\begin{equation}\label{50}
\rho\mathbf{u}\in L^\infty_tH^1_x.
\end{equation}
On the other hand, we deduce from \eqref{2.6}, \eqref{1.9}, and H{\"o}lder's inequality that
\begin{equation*}
\begin{split}
(\rho\mathbf{u})_t & =\rho_t\mathbf{u}+\rho\mathbf{u}_t \\
& =-\nabla\rho|\mathbf{u}|^2-\rho\mathbf{u}\divv\mathbf{u}
+\rho\mathbf{u}_t \\
& =-\nabla\rho\bar{x}^{a}\cdot|\mathbf{u}|^2\bar{x}^{-a}
-\rho\bar{x}^{a}\mathbf{u}\bar{x}^{-a}\divv\mathbf{u}
+\sqrt{\rho}\cdot\sqrt{\rho}\mathbf{u}_t \\
& =\left(-\nabla\rho\bar{x}^{a}\cdot(|\mathbf{u}|\bar{x}^{-\frac{a}{2}})^2
-\rho\bar{x}^{a}\mathbf{u}\bar{x}^{-a}\divv\mathbf{u}
+\sqrt{\rho}\cdot\sqrt{\rho}\mathbf{u}_t\right)\in L_{t,x}^2,
\end{split}
\end{equation*}
which together with \eqref{50} and Aubin-Lions lemma gives the continuity of $\rho\mathbf{u}$ at the initial time.
\end{remark}

\begin{remark}
It should be noted that the initial density can have compact support, due to \cite{YZ2014}, it seems to be a challenge to show the global existence of strong solutions to the problem \eqref{mhd}--\eqref{n4}.
\end{remark}

We now make some comments on the key ingredients of the analysis in this paper. One of the main difficulties in studying
such MHD system lies in the non-resistivity of the magnetic equations.
Moreover, for the whole two-dimensional space, when the far field density allows vacuum, it seems difficult to bound the $L^p(\mathbb{R}^2)$-norm of $\mathbf{u}$ just in terms of
$\|\sqrt{\rho}\mathbf{u}\|_{L^2(\mathbb{R}^2)}$ and  $\|\nabla\mathbf{u}\|_{L^2(\mathbb{R}^2)}$.
Furthermore, as mentioned in many papers (see \cite{LH2015,CZ2017} for example), the strong coupling between the velocity field and the magnetic field, such as $\||\mathbf{u}| |\mathbf{b}|\|$ and $\||\mathbf{u}||\nabla\mathbf{b}|\|$, will bring out some new difficulties.

In order to overcome these difficulties stated above, some new ideas and observations are needed. First of all, motivated by \cite{LL2014} (see also \cite{LS2015}), it is enough to bound the $L^p$-norm of the momentum $\rho\mathbf{u}$ instead of just the velocity $\mathbf{u}$. More precisely, using a Hardy-type inequality (see \eqref{2.50} and \eqref{2.6}) which is originally due to Lions \cite{L1996} with some careful estimates on the essential support of the density (see \eqref{3.31}), we first obtain a key Hardy-type inequality (see \eqref{3.32}) to bound the $L^p$-norm of $\mathbf{u}\bar{x}^{-\eta}$ instead of just the velocity $\mathbf{u}$, and then establish a crucial inequality (see \eqref{3.33}) which is used to control the $L^p$-norm of $\rho\mathbf{u}$.
Secondly, compared with the isentropic case \cite{CZ2017}, since the pressure $P$ is not a function of the density, it seems to be difficult to get some necessary lower order estimates of the solutions. It is worth noticing
that the energy equation \eqref{mhd}$_3$ is a hyperbolic one, using the basic estimates of the material derivatives of velocity, we succeed in deriving the
desired a priori estimates of the solutions.
Thirdly, inspired by \cite{Z2017}, we deduce some spatial weighted estimates on $\mathbf{b}$ (i.e., $\bar{x}^{a}\mathbf{b}$, see \eqref{3.28}) which are crucial to control the coupled terms, such as $\||\mathbf{u}||\mathbf{b}|\|$ and $\||\mathbf{u}||\nabla\mathbf{b}|\|$. Next, we then construct approximate solutions to \eqref{mhd}, that is, for density strictly away from vacuum initially, consider a initial boundary value problem of \eqref{mhd} in any bounded ball $B_R$ with radius $R>0.$  Finally, combining  all key points mentioned above with the similar arguments as in \cite{LS2015,LL2014,CZ2017}, we derive some desired bounds on the gradients of both the solutions and the spatial weighted density, which are independent of both the radius of the balls $B_R$ and the lower bound of the initial density.

The rest of the paper is organized as follows: In Section 2, we collect some
elementary facts and inequalities which will be needed in later analysis. Sections 3 is devoted to the a priori estimates which are needed to obtain  the local existence of strong solutions. The main result
Theorem \ref{t1} is proved in Section 4.

\section{Preliminaries}
In this section, we will recall some  known facts and elementary
inequalities which will be used frequently later.
First of all, if the initial density is strictly away from vacuum, the following local existence theorem on bounded balls can be shown by similar arguments as in \cite{FY2009,K1983}.
\begin{lemma}\label{lem21}
For $R>0$ and $B_R=\{x\in\mathbb{R}^2||x|<R\}$, assume that $(\rho_0,\mathbf{u}_0,P_0,\mathbf{b}_0)$ satisfies
\begin{align}\label{2.1}
\inf\limits_{x\in B_R}\rho_0(x)>0, \
(\rho_0,P_0) \in H^3(B_R), \
(\mathbf{u}_0,\mathbf{b}_0)\in H_0^1(B_R)\cap H^3(B_R), \
\divv\mathbf{b}_0=0.
\end{align}
Then there exist a small time $T_R>0$ and a unique classical solution $(\rho, \mathbf{u}, P, \mathbf{b})$ to the following initial-boundary-value conditions
\begin{align}\label{2.2}
\begin{cases}
\rho_{t}+\divv(\rho\mathbf{u})=0,\\
(\rho\mathbf{u})_{t}+\divv(\rho\mathbf{u}\otimes\mathbf{u})
-\mu\Delta\mathbf{u}
-(\lambda+\mu)\nabla\divv\mathbf{u}+\nabla P=
\mathbf{b}\cdot\nabla\mathbf{b}-\frac12\nabla|\mathbf{b}|^2,\\
P_{t}+\divv(P\mathbf{u})
+P\divv\mathbf{u}=2\mu|\mathfrak{D}(\mathbf{u})|^2+\lambda(\divv\mathbf{u})^2,\\
\mathbf{b}_{t}-\mathbf{b}\cdot\nabla\mathbf{u}
+\mathbf{u}\cdot\nabla\mathbf{b}+\mathbf{b}\divv\mathbf{u}=\mathbf{0},\\
\divv \mathbf{b}=0,\\
(\rho, \mathbf{u}, P, \mathbf{b})(x,t=0)=(\rho_0,\mathbf{u}_0, P_0, \mathbf{b}_0), \ x\in B_R,\\
\mathbf{u}(x,t),\ \mathbf{b}(x,t)=0,\  x\in\partial B_R, \, t>0,
\end{cases}
\end{align}
on $B_R\times(0,T_R]$ such that
\begin{align}\label{2.3}
\begin{cases}
\rho, P, \mathbf{b}\in C\left([0,T_{R}];H^3\right),\\
\mathbf{u} \in C\left([0,T_{R}]; H^3\right)\cap L^{2}\left(0,T_{R};H^4\right),\\
\mathbf{u}_t\in L^\infty\left(0,T_{R}; H^{1}\right)\cap L^{2}\left(0,T_{R};H^{2}\right),
\end{cases}
\end{align}
where we denote $H^k=H^k(B_R)$ for positive integer $k$.
\end{lemma}

Next, for $\Omega\subset\mathbb{R}^2$, the following weighted $L^m$-bounds for elements of the Hilbert space $\tilde{D}^{1,2}(\Omega)\triangleq\{v\in H^1_{\rm loc}(\Omega)|\nabla v\in L^2(\Omega)\}$ can be found in \cite[Theorem B.1]{L1996}.
\begin{lemma} \label{lem22}
For $m\in [2,\infty)$ and $\theta\in (1+m/2,\infty),$ there exists a positive constant $C$ such that for either $\Omega=\mathbb{R}^2$ or $\Omega=B_R$ with $R\ge 1$ and for any $v\in \tilde{D}^{1,2}(\Omega)$,
\begin{align}\label{2.4}
\left(\int_{\Omega} \frac{|v|^m}{e+|x|^2}(\log (e+|x|^2))^{-\theta}dx  \right)^{\frac1m}\le C\|v\|_{L^2(B_1)}+C\|\nabla v\|_{L^2 (\Omega)}.
\end{align}
\end{lemma}

A useful consequence of Lemma \ref{lem22} is the following crucial weighted  bounds for elements of $\tilde{D}^{1,2}(\Omega)$, which have been proved in
\cite[Lemma 2.4]{LL2014}.

\begin{lemma}\label{lem23}
Let $\bar x$ and $\eta_0$ be as in \eqref{1.6} and $\Omega$ be as in Lemma \ref{lem22}. Assume that $\rho \in L^1(\Omega)\cap L^\infty(\Omega)$ is a non-negative function such that
\begin{equation}\label{2.5}
\int_{B_{N_1}}\rho dx\geq M_1, \quad \|\rho\|_{L^1\cap L^\infty}\leq M_2,
\end{equation}
for positive constants $M_1, M_2$, and $N_1\ge 1$ with $B_{N_1}\subset\Omega$. Then
there is a positive constant $C$ depending only on $M_1,M_2,
N_1,$ and $\eta_0$ such that
\begin{align}\label{2.50}
\|v\bar x^{-1}\|_{L^{2}(\Omega)}
\leq C \|\sqrt{\rho}v\|_{L^2(\Omega)}+C \|\nabla v\|_{L^2(\Omega)}.
\end{align}
Moreover, for $\varepsilon>0$ and $\eta>0,$ there is a positive constant $C$ depending only on $\varepsilon,\eta, M_1,M_2,
N_1, $ and $\eta_0$ such that every $v\in \tilde{D}^{1,2}(\Omega)$ satisfies
\begin{align}\label{2.6}
\|v\bar x^{-\eta}\|_{L^{(2+\varepsilon)/\tilde{\eta}}(\Omega)}
\leq C \|\sqrt{\rho}v\|_{L^2(\Omega)}+C \|\nabla v\|_{L^2(\Omega)}
\end{align}
with $\tilde{\eta}=\min\{1,\eta\}$.
\end{lemma}

Finally, the following $L^p$-bound for elliptic systems, whose proof is similar to that of \cite[Lemma 12]{CCK2004}, is a direct result of the combination of the well-known elliptic theory \cite{adn} and a standard scaling procedure.
\begin{lemma}\label{lem24}
For $p>1$ and $k\ge 0,$ there exists a positive constant $C$ depending only on $p$ and $k$ such that
\begin{equation}\label{2.70}
\|\nabla^{k+2}v\|_{L^p(B_R)}\leq C\|\Delta v\|_{W^{k,p}(B_R)},
\end{equation}
for every $v\in W^{k+2,p}(B_R)$ satisfying
\begin{equation*}
v =0\ \ \mbox{on}\ \ B_R.
\end{equation*}
\end{lemma}

\section{A priori estimates}
In this section, for $r\in [1,\infty]$ and $k\ge0$, we denote
\begin{equation*}
\int \cdot dx=\int_{B_R}\cdot dx, \quad L^r=L^r(B_R),\quad W^{k,r}=W^{k,r}(B_R),\quad H^k=W^{k,2}.
\end{equation*}
Moreover, for $R>4N_0\ge 4,$ assume that $(\rho_0,\mathbf{u}_0,P_0,\mathbf{b}_0)$ satisfies, in addition to \eqref{2.1}, that
\begin{align}\label{3.1}
1/2\le \int_{B_{N_0}}\rho_0(x)dx \le \int_{B_R }\rho_0(x)dx \le 3/2.
\end{align}
Lemma \ref{lem21} thus yields that  there exists some $T_R>0$ such that the  initial-boundary-value problem \eqref{2.2} has a unique classical solution $(\rho,\mathbf{u},P,\mathbf{b})$ on $B_R\times[0,T_R]$ satisfying \eqref{2.3}.

Let $\bar{x}, \eta_0, a$, and $q$ be as in Theorem \ref{t1}, the main goal of this section is to derive the following key a priori estimate on $\psi $ defined by
\begin{align}\label{3.2}
\psi(t) & \triangleq 1+\|\sqrt{\rho}\mathbf{u}\|_{L^2}+\|\nabla\mathbf{u}\|_{L^2}
+\|\bar{x}^a\mathbf{b}\|_{L^2\cap L^{q}}+\|\nabla\mathbf{b}\|_{L^2\cap L^{q}} \notag \\
& \quad +\|\bar{x}^{a}\rho\|_{L^{1}\cap H^{1}\cap W^{1,q}}+\|P\|_{L^{1}\cap H^{1}\cap W^{1,q}}.
\end{align}
\begin{proposition} \label{pro}
Assume that $(\rho_0,\mathbf{u}_0,P_0,\mathbf{b}_0)$ satisfies \eqref{2.1} and \eqref{3.1}. Let $(\rho,\mathbf{u},P,\mathbf{b})$ be the solution to the initial-boundary-value problem \eqref{2.2} on $B_R\times (0,T_R]$ obtained by Lemma \ref{lem21}. Then there exist positive constants $T_0$ and $M$ both depending only on $\mu,q$, $a$, $\eta_0$, $N_0,$ and $E_0$ such that
\begin{align}\label{3.3}
\sup\limits_{0\le t\le T_0}\left(\psi(t)+\|\sqrt{\rho}\dot{\mathbf{u}}\|_{L^2}^2\right)
+\int_0^{T_0}\left(\|\nabla\dot{\mathbf{u}}\|_{L^2}^{2}
+\|\nabla^2\mathbf{u}\|_{L^2\cap L^q}^{2}\right)dt\leq M,
\end{align}
where
\begin{equation*}
E_0\triangleq \|\sqrt{\rho_0}\mathbf{u}_0\|_{L^2}+\|\nabla\mathbf{u}_0\|_{L^2}
+\|\bar{x}^{a}\rho_0\|_{L^{1}\cap H^{1}\cap W^{1,q}}+\|\bar{x}^a\mathbf{b}_0\|_{L^2\cap L^{q}}+\|\nabla\mathbf{b}_0\|_{L^2\cap L^{q}}+\|P_0\|_{L^{1}\cap H^{1}\cap W^{1,q}}.
\end{equation*}
\end{proposition}

To show Proposition \ref{pro}, whose proof will be postponed to the end of this section, we begin with the following standard energy estimate for $(\rho,\mathbf{u},P,\mathbf{b})$.
\begin{lemma}\label{lem32}
Under the conditions of Proposition \ref{pro}, let $(\rho,\mathbf{u},P,\mathbf{b})$ be a
smooth solution to the initial-boundary-value problem \eqref{2.2}. Then for any $t>0$,
\begin{align}\label{3.4}
\sup_{0\le s\le t} \left(\|\sqrt{\rho}\mathbf{u}\|^2_{L^{2}}+\|\mathbf{b}\|_{L^2}^2+\|P\|_{L^1}\right)
\leq C,
\end{align}
where (and in what follows) $C$ denotes a generic positive constant depending only on $\mu,\lambda, q, a$, $\eta_0$, $N_0,$ and $E_0$.
\end{lemma}
{\it Proof.}
It follows from \eqref{mhd}$_3$ that
\begin{equation}\label{3.04}
P_t+\mathbf{u}\cdot\nabla P+2P\divv\mathbf{u}
=F\triangleq2\mu|\mathfrak{D}(\mathbf{u})|^2
+\lambda(\divv\mathbf{u})^2\geq0.
\end{equation}
Define the following particle path
\begin{align*}
\begin{cases}
\frac{d}{dt}\mathbf{X}(x,t)
=\mathbf{u}(\mathbf{X}(x,t),t),\\
\mathbf{X}(x,0)=x.
\end{cases}
\end{align*}
Thus, along particle path, we obtain from \eqref{3.04} that
\begin{align*}
\frac{d}{dt}P(\mathbf{X}(x,t),t)
=-2P\divv\mathbf{u}+F,
\end{align*}
which implies for any $t>0$,
\begin{equation*}
P(\mathbf{X}(x,t),t)=\exp\left(-2\int_{0}^t\divv\mathbf{u}ds\right)
\left[P_0+\int_{0}^t\exp\left(2\int_{0}^s\divv\mathbf{u}d\tau\right)Fds\right]\geq0.
\end{equation*}

Multiplying \eqref{mhd}$_2$ and \eqref{mhd}$_4$ by $\mathbf{u}$ and $\mathbf{b}$ respectively, then adding the two resulting equations together, and integrating over $B_R$, we obtain after integrating by parts that
\begin{align}\label{3.5}
\frac{1}{2}\frac{d}{dt}\int\left(\rho|\mathbf{u}|^2
+|\mathbf{b}|^2\right)dx
+\int\left[\mu|\nabla\mathbf{u}|^2+(\lambda+\mu)(\divv\mathbf{u})^2
\right]dx
=\int P\divv\mathbf{u}dx.
\end{align}
Integrating \eqref{mhd}$_3$ with respect to $x$ and then adding the resulting equality to \eqref{3.5} give rise to
\begin{equation*}
\frac{d}{dt}\int\left(\frac12\rho|\mathbf{u}|^2+\frac12|\mathbf{b}|^2+P\right)dx=0.
\end{equation*}
Thus, integrating with respect to $t$ leads to
\begin{equation*}
\sup_{0\leq t\leq T}\left(\|\sqrt{\rho}\mathbf{u}\|_{L^2}^2+\|\mathbf{b}\|_{L^2}^2
+\|P\|_{L^1}
\right)\leq C.
\end{equation*}

This completes the proof of Lemma \ref{lem32}.   \hfill $\Box$

\begin{lemma} \label{lem33}
Let $(\rho,\mathbf{u},P,\mathbf{b})$ be a
smooth solution to  the  initial-boundary-value problem \eqref{2.2},
then there exists a positive constant $\alpha>1$ such that for any $t>0$,
\begin{align}\label{3.6}
\sup_{0\le s\le t}\left(\|\nabla\mathbf{u}\|_{L^2}^2+\|\mathbf{b}\|_{L^4}^4\right)
+\int_{0}^{t}\|\sqrt{\rho}\dot{\mathbf{u}}\|_{L^2}^2ds\le C+C\int_{0}^{t}\psi^\alpha(s)ds.
\end{align}
\end{lemma}
{\it Proof.}
Multiplying $\eqref{2.2}_{2}$ by $\dot{\mathbf{u}}$ and integrating by parts, one has
\begin{align}\label{3.7}
\int\rho|\dot{\mathbf{u}}|^2dx
& = -\int\dot{\mathbf{u}}\cdot\nabla Pdx
+\mu\int\Delta\mathbf{u}\cdot\dot{\mathbf{u}}dx
+(\mu+\lambda)\int\nabla\divv\mathbf{u}\cdot\dot{\mathbf{u}}dx \nonumber \\
& \quad -\frac12\int\dot{\mathbf{u}}\cdot\nabla|\mathbf{b}|^2dx
+\int\mathbf{b}\cdot \nabla\mathbf{b}\cdot\dot{\mathbf{u}}dx
\triangleq \sum_{i=1}^5K_i.
\end{align}
We next estimate each term on the right-hand side of \eqref{3.7} as follows.

First, integration by parts together with \eqref{2.2}$_3$ and
Gagliardo-Nirenberg inequality lead to
\begin{align}\label{3.8}
K_1 & = \int\divv\mathbf{u}_tPdx+\int\divv(\mathbf{u}\cdot\nabla\mathbf{u})Pdx \notag \\
& = \left(\int\divv\mathbf{u}Pdx\right)_t-\int\divv\mathbf{u}P_tdx
+\int\mathbf{u}\cdot\nabla(\divv\mathbf{u})Pdx+\int\partial_iu^j\partial_ju^iPdx
\notag \\
& \leq \left(\int\divv\mathbf{u}Pdx\right)_t
-\int\divv\mathbf{u}(P_t+\mathbf{u}\cdot P)dx
+C\int P|\nabla\mathbf{u}|^2dx
\notag \\
& \leq \left(\int\divv\mathbf{u}Pdx\right)_t
-\int\divv\mathbf{u}(2\mu|\mathfrak{D}(\mathbf{u})|^2+\lambda(\divv\mathbf{u})^2
-2P\divv\mathbf{u})dx+C\int P|\nabla\mathbf{u}|^2dx
\notag \\
& \leq \left(\int\divv\mathbf{u}Pdx\right)_t
+C\int|\nabla\mathbf{u}|^3dx+C\int P|\nabla\mathbf{u}|^2dx.
\end{align}
Similarly to the proof of \cite[Lemma 3.2]{LX2014}, we have
\begin{align}\label{3.9}
K_2+K_3\leq \left(-\frac{\mu}{2}\|\nabla\mathbf{u}\|_{L^2}^2
-\frac{\lambda+\mu}{2}\|\divv\mathbf{u}\|_{L^2}^2\right)_t
+C\|\nabla\mathbf{u}\|_{L^3}^3.
\end{align}
Employing \eqref{2.2}$_4$, \eqref{3.4}, and Gagliardo-Nirenberg inequality, we get
\begin{align}\label{3.10}
K_4 & = \frac12\int|\mathbf{b}|^2\divv\mathbf{u}_tdx
+\frac12\int|\mathbf{b}|^2\divv(\mathbf{u}\cdot\nabla\mathbf{u})dx \notag \\
& = \left(\frac12\int|\mathbf{b}|^2\divv\mathbf{u}dx\right)_t
-\frac12\int|\mathbf{b}|^2(\divv\mathbf{u})^2dx
+\frac12\int|\mathbf{b}|^2\partial_iu^j\partial_ju^idx
\notag \\
& \quad -\int(\mathbf{b}\cdot\nabla\mathbf{u}-\mathbf{b}\divv\mathbf{u})
\cdot\mathbf{b}\divv\mathbf{u}dx \notag \\
& \leq \left(\frac12\int|\mathbf{b}|^2\divv\mathbf{u}dx\right)_t
+C\int|\mathbf{b}|^2|\nabla\mathbf{u}|^2dx
\notag \\
& \leq \left(\frac12\int|\mathbf{b}|^2\divv\mathbf{u}dx\right)_t
+C\|\nabla\mathbf{u}\|_{L^3}^3+C\|\mathbf{b}\|_{L^6}^6
\notag \\
& \leq \left(\frac12\int|\mathbf{b}|^2\divv\mathbf{u}dx\right)_t
+C\|\nabla\mathbf{u}\|_{L^3}^3+C\|\nabla\mathbf{b}\|_{L^2}^4.
\end{align}
Similarly, one deduces
\begin{align}\label{3.11}
K_5 \leq -\frac{d}{dt}\int\mathbf{b}\cdot\nabla\mathbf{u}\cdot\mathbf{b}dx
+C\|\nabla\mathbf{u}\|_{L^3}^3+C\|\nabla\mathbf{b}\|_{L^2}^4.
\end{align}
Substituting \eqref{3.8}--\eqref{3.11} into \eqref{3.7} gives rise to
\begin{align}\label{3.12}
& \frac{d}{dt}\left(\frac{\mu}{2}\|\nabla\mathbf{u}\|_{L^2}^2
+\frac{\lambda+\mu}{2}\|\divv\mathbf{u}\|_{L^2}^2\right)
+\|\sqrt{\rho}\dot{\mathbf{u}}\|_{L^2}^2
\notag \\
& \leq B'(t)+C\int P|\nabla\mathbf{u}|^2dx
+C\|\nabla\mathbf{u}\|_{L^3}^3+C\|\nabla\mathbf{b}\|_{L^2}^4
\notag \\
& \leq B'(t)+C\|P\|_{L^\infty}\|\nabla\mathbf{u}\|_{L^2}^2
+C\|\nabla\mathbf{u}\|_{L^2}^2\|\nabla\mathbf{u}\|_{H^1}+C\|\nabla\mathbf{b}\|_{L^2}^4
\notag \\
& \leq B'(t)+C\psi^\alpha
+C\psi^2\|\nabla^2\mathbf{u}\|_{L^2},
\end{align}
where
\begin{align}\label{3.13}
B(t)\triangleq \int P\divv\mathbf{u}dx+\frac12\int\divv\mathbf{u}|\mathbf{b}|^2dx
-\int\mathbf{b}\cdot\nabla\mathbf{u}\cdot\mathbf{b}dx.
\end{align}
Here (and in what follows) we use $\alpha>1$ to denote a genetic constant, which may be different from line to line.

To estimate the last term on the right-hand side of \eqref{3.12}, it follows from \eqref{2.70} and \eqref{2.2}$_2$ that for any $p\geq2$,
\begin{equation}\label{3.14}
\|\nabla^2\mathbf{u}\|_{L^p}\leq C\|\rho\dot{\mathbf{u}}\|_{L^p}
+C\|\nabla P\|_{L^p}+C\||\mathbf{b}||\nabla\mathbf{b}|\|_{L^p}.
\end{equation}
Choosing $p=2$ in \eqref{3.14}, then we get from \eqref{3.12} that
\begin{align}\label{3.15}
& \frac{d}{dt}\left(\frac{\mu}{2}\|\nabla\mathbf{u}\|_{L^2}^2
+\frac{\lambda+\mu}{2}\|\divv\mathbf{u}\|_{L^2}^2\right)
+\|\sqrt{\rho}\dot{\mathbf{u}}\|_{L^2}^2
\notag \\
& \leq B'(t)+C\psi^\alpha
+C\psi^2\left(\|\rho\|_{L^\infty}^{\frac12}\|\sqrt{\rho}\dot{\mathbf{u}}\|_{L^2}
+C\|\nabla P\|_{L^2}+C\|\mathbf{b}\|_{L^\infty}\|\nabla\mathbf{b}\|_{L^2}\right)
\notag \\
& \leq B'(t)+C\psi^\alpha
+\frac12\|\sqrt{\rho}\dot{\mathbf{u}}\|_{L^2}^2.
\end{align}
Integrating \eqref{3.15} over $(0,t)$ and using \eqref{3.13}, one arrives at
\begin{align}\label{3.16}
& \sup_{0\leq s\leq t}\|\nabla\mathbf{u}\|_{L^2}^2
+\int_{0}^t\|\sqrt{\rho}\dot{\mathbf{u}}\|_{L^2}^2ds \notag \\
& \leq C+C\int_{0}^t\psi^\alpha ds
+\int\divv\mathbf{u}Pdx+\frac12\int\divv\mathbf{u}|\mathbf{b}|^2dx
-\int\mathbf{b}\cdot\nabla\mathbf{u}\cdot\mathbf{b}dx \notag \\
& \leq C+C\int_{0}^t\psi^\alpha ds
+\frac12\|\nabla\mathbf{u}\|_{L^2}^2+C\|P\|_{L^2}^2+C\|\mathbf{b}\|_{L^4}^4
 \notag \\
& \leq C+C\int_{0}^t\psi^\alpha ds
+\frac12\|\nabla\mathbf{u}\|_{L^2}^2+C\|P\|_{L^2}^2,
\end{align}
where in the last inequality we have used the following
\begin{align}\label{3.17}
\sup_{0\leq s\leq t}\|\mathbf{b}\|_{L^4}^4
\leq C+C\int_{0}^t\psi^\alpha ds.
\end{align}
Indeed, multiplying \eqref{mhd}$_4$ by $4|\mathbf{b}|^2\mathbf{b}$ and integrating the resulting equation over $B_R$, we derive that
\begin{align*}
\frac{d}{dt}\|\mathbf{b}\|_{L^4}^4
& \leq C\int|\nabla\mathbf{u}||\mathbf{b}|^4dx \\
& \leq C\|\nabla\mathbf{u}\|_{L^2}\|\mathbf{b}\|_{L^8}^4 \\
& \leq C\|\nabla\mathbf{u}\|_{L^2}\|\mathbf{b}\|_{L^2}\|\nabla\mathbf{b}\|_{L^2}^3 \\
& \leq C\psi^\alpha.
\end{align*}
Integrating the above inequality over $(0,t)$ yields \eqref{3.17}.

Finally, in order to estimate the last term on the right-hand side of \eqref{3.16}, multiplying \eqref{2.2}$_3$ by $2P$ and noting that
\begin{align*}
2\int\divv(P\mathbf{u})Pdx
& = 2\int\divv(P^2\mathbf{u})dx-2\int P\mathbf{u}\cdot\nabla Pdx \\
& = -\int \mathbf{u}\cdot\nabla P^2dx
= \int P^2\divv\mathbf{u}dx,
\end{align*}
we obtain that
\begin{align*}
\frac{d}{dt}\|P\|_{L^2}^2
& \leq C\int P^2|\divv\mathbf{u}|dx+C\int P|\nabla\mathbf{u}|^2dx \notag \\
& \leq C\|P\|_{L^\infty}\|P\|_{L^2}\|\nabla\mathbf{u}\|_{L^2}
+C\|P\|_{L^\infty}\|\nabla\mathbf{u}\|_{L^2}^2 \notag \\
& \leq C\psi^\alpha,
\end{align*}
which leads to
\begin{equation}\label{3.18}
\sup_{0\leq s\leq t}\|P\|_{L^2}^2
\leq C+C\int_{0}^t\psi^\alpha ds.
\end{equation}
So the desired \eqref{3.6} follows from \eqref{3.16}, \eqref{3.17}, and \eqref{3.18}. The proof of Lemma \ref{lem33} is finished.  \hfill $\Box$

\begin{lemma} \label{lem34}
Let $\alpha$ be as in Lemma \ref{lem33}. Then for any $t>0$,
\begin{align}\label{3.19}
\sup_{0\le s\le t}\|\sqrt{\rho}\dot{\mathbf{u}}\|_{L^2}^2
+\int_{0}^{t}\|\nabla\dot{\mathbf{u}}\|_{L^2}^2ds
\leq C\exp\left\{C\int_{0}^{t} \psi^\alpha ds\right\}.
\end{align}
\end{lemma}
{\it Proof.}
Operating $\partial_t+\divv(\mathbf{u}\cdot)$ to the $j$-th component of \eqref{2.2}$_2$
and multiplying the resulting equation by $\dot{u}^j$, one gets by some calculations that
\begin{align}\label{3.20}
\frac12\frac{d}{dt}\int \rho|\dot{\mathbf{u}}|^2\mbox{d}x
& =\mu\int\dot{u}^j(\partial_t\Delta u^j+\divv(\mathbf{u}\Delta u^j))dx
+(\lambda+\mu)\int \dot{u}^j (\partial_t\partial_j(\divv\mathbf{u})+\divv(\mathbf{u}\partial_j(\divv\mathbf{u})))dx  \nonumber\\
& \quad -\int \dot{u}^j(\partial_jP_t+\divv(\mathbf{u}\partial_jP)) dx
-\frac12\int \dot{u}^j (\partial_t\partial_j|\mathbf{b}|^2+\divv(\mathbf{u}\partial_j|\mathbf{b}|^2)) dx\nonumber\\
& \quad +\int \dot{u}^j (\partial_t(\mathbf{b}\cdot\nabla b^j)+\divv(\mathbf{u}(\mathbf{b}\nabla b^j))) dx
\triangleq\sum_{i=1}^{5}J_i,
\end{align}
where $J_i$ can be bounded as follows.

Integration by parts leads to
\begin{align}\label{3.21}
J_1 & = - \mu\int(\partial_i\dot{u}^j\partial_t\partial_i u^j+\Delta u^j\mathbf{u}\cdot\nabla\dot{u}^j)dx \notag \\
& = - \mu\int(|\nabla\dot{\mathbf{u}}|^2
-\partial_i\dot{u}^ju^k\partial_k\partial_i u^j
-\partial_i\dot{u}^j\partial_iu^k\partial_k u^j
+\Delta u^j\mathbf{u}\cdot\nabla\dot{u}^j)dx
\notag \\
& = - \mu\int(|\nabla\dot{\mathbf{u}}|^2
+\partial_i\dot{u}^j\partial_k u^k\partial_i u^j
-\partial_i\dot{u}^j\partial_iu^k\partial_k u^j
-\partial_i u^j\partial_iu^k\partial_k\dot{u}^j)dx
\notag \\
& \leq -\frac{3\mu}{4}\|\nabla\dot{\mathbf{u}}\|_{L^2}^2+C\|\nabla\mathbf{u}\|_{L^4}^4.
\end{align}
Similarly, one has
\begin{align}\label{3.22}
J_2 \leq -\frac{\lambda+\mu}{2}\|\divv\dot{\mathbf{u}}\|_{L^2}^2
+C\|\nabla\mathbf{u}\|_{L^4}^4.
\end{align}
It follows from integration by parts, \eqref{mhd}$_3$, and \eqref{3.6} that
\begin{align}\label{3.23}
J_3 & = \int(\partial_j\dot{u}^jP_t+\partial_j P\mathbf{u}\cdot\nabla\dot{u}^j)dx \notag \\
& = \int\partial_j\dot{u}^j\left(
2\mu|\mathfrak{D}(\mathbf{u})|^2+\lambda(\divv\mathbf{u})^2
-\divv(P\mathbf{u})
-P\divv\mathbf{u}\right)dx-\int P\partial_j(\mathbf{u}\cdot\nabla\dot{u}^j)dx\notag \\
& = \int\partial_j\dot{u}^j\left(
2\mu|\mathfrak{D}(\mathbf{u})|^2+\lambda(\divv\mathbf{u})^2
-\divv(P\mathbf{u})
-P\divv\mathbf{u}\right)dx \notag \\
& \quad
-\int P\partial_j\mathbf{u}\cdot\nabla\dot{u}^jdx
+\int\partial_j\dot{u}^j \divv(P\mathbf{u})dx
\notag \\
& = \int\partial_j\dot{u}^j\left(
2\mu|\mathfrak{D}(\mathbf{u})|^2+\lambda(\divv\mathbf{u})^2
-P\divv\mathbf{u}\right)dx-\int P\partial_j\mathbf{u}\cdot\nabla\dot{u}^jdx
\notag \\
& \leq C\int |\nabla\dot{\mathbf{u}}|(|\nabla\mathbf{u}|^2+P|\nabla\mathbf{u}|)dx
\notag \\
& \leq \frac{\mu}{4}\|\nabla\dot{\mathbf{u}}\|_{L^2}^2
+C\|\nabla\mathbf{u}\|_{L^4}^4+C\|P\|_{L^4}^4.
\end{align}
From \eqref{2.2}$_4$ and \eqref{2.2}$_5$, we arrive at
\begin{align}\label{3.24}
J_4 & = \int\partial_j\dot{u}^j\mathbf{b}\cdot\mathbf{b}_tdx +\frac12\int\mathbf{u}\cdot\nabla\dot{u}^j\partial_j|\mathbf{b}|^2dx \notag \\
& = \int\partial_j\dot{u}^j\mathbf{b}\cdot(\mathbf{b}\cdot\nabla\mathbf{u}
-\mathbf{b}\divv\mathbf{u})dx +\frac12\int\partial_j\dot{u}^j\divv\mathbf{u}|\mathbf{b}|^2dx
-\frac12\int\partial_j\mathbf{u}\cdot\nabla\dot{u}^j|\mathbf{b}|^2dx \notag \\
& \leq C\int |\nabla\dot{\mathbf{u}}||\nabla\mathbf{u}||\mathbf{b}|^2dx
\notag \\
& \leq \frac{\mu}{8}\|\nabla\dot{\mathbf{u}}\|_{L^2}^2
+C\|\nabla\mathbf{u}\|_{L^4}^4+C\||\mathbf{b}|^2\|_{L^4}^4.
\end{align}
Similarly to $J_4$, we infer that
\begin{align}\label{3.25}
J_5 \leq \frac{\mu}{8}\|\nabla\dot{\mathbf{u}}\|_{L^2}^2
+C\|\nabla\mathbf{u}\|_{L^4}^4+C\||\mathbf{b}|^2\|_{L^4}^4.
\end{align}
Inserting \eqref{3.21}--\eqref{3.25} into \eqref{3.20} yields
\begin{align}\label{3.26}
\frac{d}{dt}\|\sqrt{\rho}\dot{\mathbf{u}}\|_{L^2}^2
+\mu\|\nabla\dot{\mathbf{u}}\|_{L^2}^2
\leq C\|\nabla\mathbf{u}\|_{L^4}^4+C\|P\|_{L^4}^4
+C\||\mathbf{b}|^2\|_{L^4}^4.
\end{align}
Thus we obtain from Gagliardo-Nirenberg inequality, \eqref{3.14}, and \eqref{3.2} that
\begin{align}\label{3.27}
\frac{d}{dt}\|\sqrt{\rho}\dot{\mathbf{u}}\|_{L^2}^2
+\mu\|\nabla\dot{\mathbf{u}}\|_{L^2}^2
& \leq C\|\nabla\mathbf{u}\|_{L^4}^4+C\|P\|_{L^4}^4
+C\||\mathbf{b}|^2\|_{L^4}^4 \notag \\
& \leq C\|\nabla\mathbf{u}\|_{L^2}^2\|\nabla\mathbf{u}\|_{H^1}^2+C\psi^\alpha
\notag \\
& \leq C\psi^\alpha\|\nabla^2\mathbf{u}\|_{L^2}^2+C\psi^\alpha
\notag \\
& \leq C\psi^\alpha\|\sqrt{\rho}\dot{\mathbf{u}}\|_{L^2}^2+C\psi^\alpha,
\end{align}
which combined with Gronwall's inequality and \eqref{1.8} leads to the desired \eqref{3.19} and finishes the proof of Lemma \ref{lem34}.
\hfill $\Box$

\begin{lemma}\label{lem35}
Let $\alpha$ be as in Lemma \ref{lem33}. Then there exists a positive constant $T_1=T_1(N_0,E_0)$ such that for all $t\in(0, T_1]$,
\begin{align}\label{3.28}
\sup\limits_{0\leq s\leq t} \left(\|\rho \bar x^a\|_{L^1\cap H^1\cap W^{1,q}}+\|\mathbf{b}\bar{x}^a\|_{L^2\cap L^{q}}
+\|\nabla\mathbf{b}\|_{L^2\cap L^{q}}\right)
\leq \exp\left\{C\exp\left\{C\int_{0}^{t} \psi^{\alpha} ds\right\}\right\}.
\end{align}
\end{lemma}
{\it Proof.}
1. For $N>1$, let $\varphi_N\in C^\infty_0(B_N)$ satisfy
\begin{align}\label{3.29}
0\leq\varphi_N \le 1,\ \varphi_N(x)=1,\ \mbox{if}\ |x|\le N/2, \quad |\nabla\varphi_N|\le C N^{-1}.
\end{align}
It follows from \eqref{2.2}$_1$ and  \eqref{3.1} that
\begin{align}\label{3.30}
\frac{d}{dt}\int \rho\varphi_{2N_0}dx
& =\int \rho\mathbf{u}\cdot\nabla\varphi_{2N_0}dx  \notag \\
& \ge - C N_0^{-1}\left(\int\rho dx\right)^{1/2}
\left(\int\rho|\mathbf{u}|^2dx\right)^{1/2}\ge -\tilde{C}(E_0),
\end{align}
where in the last inequality we have used
\begin{equation*}
\int\rho dx=\int\rho_0dx
\end{equation*}
due to \eqref{2.2}$_1$.
Integrating \eqref{3.30} gives  rise to
\begin{align}\label{3.31}
\inf\limits_{0\le t\le T_1}\int_{B_{2N_0}}  \rho dx
& \ge \inf\limits_{0\le t\le T_1}\int \rho \varphi_{2N_0} dx\ge \int \rho_0\varphi_{2N_0} dx-\tilde{C}T_1 \ge 1/4.
\end{align}
Here, $T_1\triangleq\min\{1, (4\tilde{C})^{-1}\}$.
From now on, we will always assume that $t\le T_1.$
The combination of \eqref{3.31}, \eqref{3.4}, and \eqref{2.6} implies that for $\varepsilon>0$ and $\eta>0$, every $v\in \tilde{D}^{1,2}(B_R)$ satisfies
\begin{align}\label{3.32}
\|v\bar x^{-\eta}\|_{L^{(2+\varepsilon)/\tilde{\eta}}}^2
& \le C(\varepsilon,\eta) \|\sqrt{\rho}v\|_{L^2}^2+C(\varepsilon,\eta)
\|\nabla v\|_{L^2 }^2,
\end{align}
with $\tilde{\eta}=\min\{1,\eta\}$.
In particular, we have
\begin{align}\label{3.33}
\|\rho^\eta\mathbf{u}\|_{L^{(2+\varepsilon)/\tilde{\eta}}}
+\|\mathbf{u}\bar x^{-\eta}\|_{L^{(2+\varepsilon)/\tilde{\eta}}}\leq C(\varepsilon,\eta)\psi^{1+\eta}.
\end{align}

2. Multiplying \eqref{2.2}$_1$ by $\bar{x}^{a}$ and integrating by parts yield
\begin{align*}
\frac{d}{dt}\int\rho\bar{x}^{a}dx
& \leq C\int\rho|\mathbf{u} |\bar{x}^{a-1}\log^{1+\eta_0}(e+|x|^2)dx \notag  \\
& = C\int\rho\bar{x}^{a-1+\frac{8}{8+a}}|\mathbf{u}| \bar{x}^{-\frac{4}{8+a}}\left(\bar{x}^{-\frac{4}{8+a}}\log^{1+\eta_0}(e+|x|^2)\right)dx \notag  \\
 &  \leq C\|\rho\bar{x}^{a-1+\frac{8}{8+a}}\|_{L^{\frac{8+a}{7+a}}}\|\mathbf{u} \bar{x}^{-\frac{4}{8+a}}\|_{L^{8+a}} \notag \\
 &  \leq C\psi^\alpha
\end{align*}
owing to \eqref{3.33}. This leads to
\begin{align}\label{3.34}
\sup\limits_{0\leq s\leq t}\|\rho \bar x^a\|_{L^1}
\leq C\exp\left\{C\int_{0}^{t} \psi^{\alpha} ds\right\}.
\end{align}

3. It follows from Sobolev's inequality and \eqref{3.33} that
for $0<\delta<1,$
\begin{align}\label{3.35}
\norm[L^{\infty}]{\mathbf{u}\bar x^{-\delta}}
& \leq C(\delta)\left(\norm[L^{\frac4\delta}]{\mathbf{u}\bar x^{-\delta}}+\norm[L^{3}]{\nabla (\mathbf{u}\bar x^{-\delta})}\right) \notag \\
& \leq C(\delta)\left(\norm[L^{\frac4\delta}]{\mathbf{u}\bar x^{-\delta}}
+\norm[L^{3}]{\nabla \mathbf{u}}+\norm[L^{\frac4\delta}]{\mathbf{u}\bar x^{-\delta}}
\|\bar{x}^{-1}\nabla\bar{x}\|_{L^{\frac{12}{4-3\delta}}} \right) \notag \\
& \leq C(\delta) \left(\psi^\alpha+\norm[L^2]{\nabla^2\mathbf{u}} \right).
\end{align}
Now, One derives from \eqref{2.2}$_1$ that $\rho\bar{x}^a$ satisfies
\begin{align}\label{3.36}
\partial_{t}(\rho\bar{x}^a)+\mathbf{u}\cdot\nabla(\rho\bar{x}^a)
-a\rho\bar{x}^{a}\mathbf{u}\cdot\nabla\log\bar{x}+\rho\bar{x}^a\divv\mathbf{u}=0,
\end{align}
which along with \eqref{3.35} gives that
for any $r\in[2,q]$,
\begin{align}\label{3.37}
\frac{d}{dt}\|\nabla(\rho\bar{x}^a)\|_{L^r} \leq
& C\left(1+\|\nabla\mathbf{u}\|_{L^\infty}
+\|\mathbf{u}\cdot\nabla\log\bar{x}\|_{L^\infty}\right)
\|\nabla(\rho\bar{x}^a)\|_{L^r} \notag \\
& +C\|\rho\bar{x}^a\|_{L^\infty}\left(\||\nabla\mathbf{u} ||\nabla\log\bar{x}|\|_{L^r}+\||\mathbf{u}||\nabla^{2}\log\bar{x}|\|_{L^r}\right) \notag \\ \leq & C\left(\psi^\alpha+\|\nabla^2\mathbf{u}\|_{L^2\cap L^q}\right)\|\nabla(\rho\bar{x}^a)\|_{L^r} \notag \\
& +C\|\rho\bar{x}^a\|_{L^\infty}\left(\|\nabla\mathbf{u}\|_{L^r}
+\|\mathbf{u}\bar{x}^{-\frac{2}{5}}\|_{L^{4r}}\|\bar{x}^{-\frac{3}{2}}\|_{L^{\frac{4r}{3}}}
\right) \notag \\
\leq & C\left(\psi^\alpha+\|\nabla^2\mathbf{u}\|_{L^2\cap L^q}\right)
\left(1+\|\nabla(\rho\bar{x}^a)\|_{L^r}+\|\nabla(\rho\bar{x}^a)\|_{L^q}\right).
\end{align}
We claim that we have obtained the following estimate
\begin{align}\label{3.38}
\int_0^t\|\nabla^2\mathbf{u}\|_{L^2\cap L^q}^2ds
\leq C\exp\left\{C\int_{0}^{t} \psi^{\alpha} ds\right\},
\end{align}
whose proof will be given later.
Then we derive from \eqref{3.37}, \eqref{3.38}, and Gronwall's inequality that
\begin{align}\label{3.39}
\sup\limits_{0\leq s\leq t}\|\rho \bar x^a\|_{H^1\cap W^{1,q}}
\leq \exp\left\{C\exp\left\{C\int_{0}^{t} \psi^{\alpha} ds\right\}\right\}.
\end{align}

4. It follows from \eqref{2.2}$_4$ that $\mathbf{b}\bar{x}^a$ satisfies
\begin{align}\label{3.40}
\partial_{t}(\mathbf{b}\bar{x}^a)+\mathbf{u}\cdot\nabla(\mathbf{b}\bar{x}^a)
-a\mathbf{u}\cdot\nabla\log\bar{x}\cdot\mathbf{b}\bar{x}^{a}
+\mathbf{b}\bar{x}^a\divv\mathbf{u}=\mathbf{b}\bar{x}^a\cdot\nabla\mathbf{u},
\end{align}
which combined with \eqref{3.35} and Gagliardo-Nirenberg inequality implies that for any $r\in[2,q]$,
\begin{align}\label{3.41}
\frac{d}{dt}\|\mathbf{b}\bar{x}^a\|_{L^r}
\leq & C\left(\|\nabla\mathbf{u}\|_{L^\infty}
+\|\mathbf{u}\cdot\nabla\log\bar{x}\|_{L^\infty}\right)
\|\mathbf{b}\bar{x}^a\|_{L^r} \notag \\
\leq & C\left(\psi^\alpha+\|\nabla^2\mathbf{u}\|_{L^2\cap L^q}\right)\|\mathbf{b}\bar{x}^a\|_{L^r}.
\end{align}
Moreover, operating $\nabla$ to \eqref{2.2}$_4$ and then multiplying
$|\nabla \mathbf{b}|^{r-2}\nabla \mathbf{b}$ for $r\in[2,q]$ gives that
\begin{align}\label{3.42}
\frac{d}{dt}\|\nabla\mathbf{b}\|_{L^r} \leq
& C\|\nabla\mathbf{u}\|_{L^\infty}
\|\nabla\mathbf{b}\|_{L^r}+C\|\mathbf{b}\|_{L^\infty}\|\nabla^2\mathbf{u}\|_{L^r} \notag \\ \leq & C\left(\psi^\alpha+\|\nabla^2\mathbf{u}\|_{L^2\cap L^q}\right)
\|\nabla\mathbf{b}\|_{L^r}+C\psi^\alpha
+\|\nabla^2\mathbf{u}\|_{L^2\cap L^q}^{2} \notag \\
\leq & C\left(\psi^\alpha+\|\nabla^2\mathbf{u}\|_{L^2\cap L^q}^{2}\right)
\left(1+\|\nabla\mathbf{b}\|_{L^r}\right).
\end{align}
Combining \eqref{3.41} and \eqref{3.42}, one has
\begin{align*}
\frac{d}{dt}\left(\|\mathbf{b}\bar{x}^a\|_{L^r}
+\|\nabla\mathbf{b}\|_{L^r}\right)
\leq C\left(\psi^\alpha+\|\nabla^2\mathbf{u}\|_{L^2\cap L^q}^{\frac{q+1}{q}}\right)
\left(1+\|\mathbf{b}\bar{x}^a\|_{L^r}+\|\nabla\mathbf{b}\|_{L^r}\right),
\end{align*}
which along with Gronwall's inequality and \eqref{3.38} leads to
\begin{align}\label{3.43}
\sup\limits_{0\leq s\leq t}
\left(\|\mathbf{b}\bar{x}^a\|_{L^2\cap L^q}
+\|\nabla\mathbf{b}\|_{L^2\cap L^q}\right)
\leq \exp\left\{C\exp\left\{C\int_{0}^{t} \psi^{\alpha} ds\right\}\right\}.
\end{align}
Hence the desired \eqref{3.28} follows from \eqref{3.34}, \eqref{3.39}, and \eqref{3.43}.

5. To finish the proof of Lemma \ref{lem35}, it remains to show \eqref{3.38}.
Indeed, for any $r\in[2,q]$, one infers from \eqref{3.14}, \eqref{3.33}, and Gagliardo-Nirenberg inequality that
\begin{align}\label{3.44}
\|\nabla^2\mathbf{u}\|_{L^r}^2
& \leq C\|\rho\dot{\mathbf{u}}\|_{L^r}^2
+C\|\nabla P\|_{L^r}^2+C\||\mathbf{b}||\nabla\mathbf{b}|\|_{L^r}^2 \notag \\
& \leq C\|\rho\dot{\mathbf{u}}\|_{L^2}^{\frac{2(r^2-2r)}{r^2-2}}
\|\rho\dot{\mathbf{u}}\|_{L^{r^2}}^{\frac{4(r-1)}{r^2-2}}
+C\psi^\alpha+C\|\mathbf{b}\|_{L^\infty}^2\|\nabla\mathbf{b}\|_{L^r}^2 \notag \\
& \leq C\psi^\alpha\left(\|\sqrt{\rho}\dot{\mathbf{u}}\|_{L^2}^2\right)^{\frac{2(r-1)}{r^2-2}}
\left(\|\sqrt{\rho}\dot{\mathbf{u}}\|_{L^2}^2
+(1+\|\rho\|_{L^\infty})^2\|\nabla\dot{\mathbf{u}}\|_{L^2}^2\right)^{\frac{r^2-2r}{r^2-2}}
+C\psi^\alpha \notag \\
& \leq C\psi^\alpha\|\sqrt{\rho}\dot{\mathbf{u}}\|_{L^2}^2+C\|\nabla\dot{\mathbf{u}}\|_{L^2}^2
+C\psi^\alpha,
\end{align}
which along with \eqref{3.6} and \eqref{3.19} implies \eqref{3.38}.
\hfill $\Box$

\begin{lemma}\label{lem36}
Let $\alpha$ and $T_1$ be as in Lemmas \ref{lem33} and \ref{lem35}. Then for all $t\in (0, T_1]$,
\begin{align}\label{3.45}
\sup\limits_{0\leq s\leq t}\|P\|_{H^1\cap W^{1,q}}
\leq \exp\left\{C\exp\left\{C\int_{0}^{t} \psi^{\alpha} ds\right\}\right\}.
\end{align}
\end{lemma}
{\it Proof.}
Operating $\nabla$ to \eqref{2.2}$_3$, one derives
\begin{align}\label{3.46}
(\nabla P)_{t}+\nabla(\mathbf{u}\cdot\nabla P)
+2\nabla(P\divv\mathbf{u})=\nabla\left(2\mu|\mathfrak{D}(\mathbf{u})|^2
+\lambda(\divv\mathbf{u})^2\right),
\end{align}
which multiplied by $|\nabla P|^{r-2}\nabla P$ for $r\in[2,q]$ gives
\begin{align}\label{3.47}
\frac1r\frac{d}{dt}\int|\nabla P|^rdx
& = -\int\nabla(\mathbf{u}\cdot\nabla P)\cdot\nabla P|\nabla P|^{r-2}dx
-2\int\nabla(P\divv\mathbf{u})\cdot\nabla P|\nabla P|^{r-2}dx \notag \\
& \quad +\int\nabla\left(2\mu|\mathfrak{D}(\mathbf{u})|^2
+\lambda(\divv\mathbf{u})^2\right)\cdot\nabla P|\nabla P|^{r-2}dx\triangleq \sum_{i=1}^3\bar{K_i},
\end{align}
where $\bar{K_i}$ can be bounded as follows
\begin{align*}
\bar{K_1} & = -\int\nabla\mathbf{u}\cdot\nabla P\cdot\nabla P|\nabla P|^{r-2}dx
-\int\mathbf{u}\cdot\nabla(\nabla P)\cdot\nabla P|\nabla P|^{r-2}dx
\\ & = -\int\nabla\mathbf{u}\cdot\nabla P\cdot\nabla P|\nabla P|^{r-2}dx
-\frac1r\int\mathbf{u}\cdot\nabla(|\nabla P|^r)dx \\
& \leq C\int|\nabla\mathbf{u}||\nabla P|^{r}dx\leq C\|\nabla\mathbf{u}\|_{L^\infty}\|\nabla P\|_{L^r}^{r}, \\
\bar{K_2} & \leq C\int|\nabla\mathbf{u}||\nabla P|^{r}dx
+C\int P|\nabla^2\mathbf{u}||\nabla P|^{r-1}dx \\
& \leq C\|\nabla\mathbf{u}\|_{L^\infty}\|\nabla P\|_{L^r}^{r}
+C\|P\|_{L^\infty}\|\nabla^2\mathbf{u}\|_{L^r}\|\nabla P\|_{L^r}^{r-1}, \\
\bar{K_3} & \leq C\int|\nabla\mathbf{u}||\nabla^2\mathbf{u}||\nabla P|^{r-1}dx \\
& \leq C\|\nabla\mathbf{u}\|_{L^\infty}\|\nabla^2\mathbf{u}\|_{L^r}\|\nabla P\|_{L^r}^{r-1}.
\end{align*}
Thus we get after inserting the above estimates into \eqref{3.47} that
\begin{align}\label{3.48}
\frac{d}{dt}\|\nabla P\|_{L^2\cap L^q}
& \leq C\left(\psi^\alpha+\|\nabla^2\mathbf{u}\|_{L^2\cap L^q}\right)
\left(1+\|\nabla P\|_{L^2\cap L^q}\right)
+C\|\nabla\mathbf{u}\|_{L^\infty}\|\nabla^2\mathbf{u}\|_{L^2\cap L^q}.
\end{align}
One deduces directly from Sobolev's inequality and \eqref{3.44} that
\begin{align*}
\int\|\nabla\mathbf{u}\|_{L^\infty}\|\nabla^2\mathbf{u}\|_{L^2\cap L^q}ds
& \leq C\int\left(\psi^\alpha+\|\nabla^2\mathbf{u}\|_{L^2}^2
+\|\nabla^2\mathbf{u}\|_{L^q}^2\right)ds \\
& \leq C\exp\left\{C\int_{0}^{t} \psi^{\alpha} ds\right\},
\end{align*}
which together with \eqref{3.48} and Gronwall's inequality gives the desired \eqref{3.45}.
This completes the proof of Lemma \ref{lem36}.
\hfill $\Box$

Now, Proposition \ref{pro} is a direct consequence of  Lemmas \ref{lem32}--\ref{lem36}.

\textbf{Proof of Proposition \ref{pro}}.
It follows from \eqref{3.4}, \eqref{3.6}, \eqref{3.28}, and \eqref{3.45}  that
\begin{align*}
\psi(t)\leq \exp\left\{C\exp\left\{C\int_{0}^{t} \psi^{\alpha} ds\right\}\right\}.
\end{align*}
Standard arguments yield  that for $M\triangleq e^{Ce}$ and $T_0\triangleq \min\{T_1,(CM^\alpha)^{-1}\}$,
\begin{equation*}
\sup\limits_{0\le t\le T_0}\psi(t)\le M,
\end{equation*}
which together with \eqref{3.19} and \eqref{3.38} gives \eqref{3.3}. The proof of Proposition \ref{pro} is completed.   \hfill $\Box$

\section{Proof of Theorem \ref{t1}}

With the a priori estimates in Section 3 at hand, it is a position to show Theorem 1.1.

\textbf{Proof of Theorem \ref{t1}.}
Let $(\rho_{0},\mathbf{u}_{0},P_0,\mathbf{b}_{0})$ be as in Theorem \ref{t1}. Without loss of generality, we assume that the initial density $\rho_0$ satisfies
\begin{equation*}
\int_{\mathbb{R}^2} \rho_0dx=1,
\end{equation*}
which implies that there exists a positive constant $N_0$ such that
\begin{equation}\label{4.1}
\int_{B_{N_0}} \rho_0 dx\ge \frac34\int_{\mathbb{R}^2}\rho_0dx=\frac34.
\end{equation}
We construct
$\rho_{0}^{R}=\hat\rho_{0}^{R}+R^{-1}e^{-|x|^2}$, where $0\le\hat\rho_{0}^{R}\in  C^\infty_0(\mathbb{R}^2)$ satisfies
\begin{align}\label{4.2}
\begin{cases}
\int_{B_{N_0}}\hat\rho^R_0dx\ge 1/2,\\
\bar x^a \hat\rho_{0}^{R}\rightarrow \bar x^a \rho_{0}\quad \text{in}\,\, L^1(\mathbb{R}^2)\cap H^{1}(\mathbb{R}^2)\cap W^{1,q}(\mathbb{R}^2) ,\,\,{\rm as}\,\,R\rightarrow\infty.
\end{cases}
\end{align}
Similarly, we can also choose
$P_{0}^{R}\in  C^\infty_0(\mathbb{R}^2)$ such that
\begin{align}\label{4.3}
P_{0}^{R}\rightarrow P_{0}\quad \text{in}\,\, L^1(\mathbb{R}^2)\cap H^{1}(\mathbb{R}^2)\cap W^{1,q}(\mathbb{R}^2) ,\,\,{\rm as}\,\,R\rightarrow\infty.
\end{align}
Then we choose $\mathbf{b}_0^R\in \{\mathbf{w}\in C^\infty_0(B_R)~|\divv\mathbf{w}=0\}$ satisfying
\begin{align}\label{4.4}
\mathbf{b}_0^R\bar x^a \rightarrow  \mathbf{b}_0\bar x^a \quad
{\rm in}\,\, L^2(\mathbb{R}^2)\cap L^{q}(\mathbb{R}^2),\quad
\nabla\mathbf{b}_0^R \rightarrow  \nabla\mathbf{b}_0 \quad
{\rm in}\,\, L^2(\mathbb{R}^2)\cap L^{q}(\mathbb{R}^2),
\quad{\rm as}\,\,R\rightarrow\infty.
\end{align}

We consider the unique smooth solution $\mathbf{u}_0^R$ of the following elliptic problem:
\begin{align}\label{4.5}
\begin{cases}
-\mu\Delta\mathbf{u}_{0}^{R}-(\mu+\lambda)\nabla\divv\mathbf{u}_{0}^{R}+\nabla P^{R}_0 =\mathbf{b}_0^R\cdot\nabla\mathbf{b}_0^R-\frac12\nabla|\mathbf{b}_0^R|^2
+\sqrt{\rho_{0}^{R}}\mathbf{h}^R-\rho_0^R\mathbf{u}_0^R,\\
\mathbf{u}_{0}^{R}|_{\partial B_{R}}=\mathbf{0},
\end{cases}
\end{align}
where $\mathbf{h}^R=(\sqrt{\rho_0}\mathbf{u}_0)*j_{1/R}$ with $j_\delta$ being the standard mollifying kernel of width $\delta$.
Extending $\mathbf{u}_{0}^{R} $ to $\mathbb{R}^2$ by defining $0$ outside $B_{R}$ and denoting it by $\tilde{\mathbf{u}}_{0}^{R}$, we claim that
\begin{align}\label{4.6}
\lim\limits_{R\rightarrow \infty}
\left(\|\nabla(\tilde{\mathbf{u}}_{0}^{R}-\mathbf{u}_0)\|_{L^2(\mathbb{R}^2)}
+\|\sqrt{\rho_0^R}\tilde{\mathbf{u}}_{0}^{R}
-\sqrt{\rho_0}\mathbf{u}_0\|_{L^2(\mathbb{R}^2)}\right)=0.
\end{align}
In fact, it is easy to find that $\tilde{\mathbf{u}}_0^R$ is also a solution of \eqref{4.5} in $\mathbb{R}^2$. Multiplying \eqref{4.5} by $\tilde{\mathbf{u}}_0^R$ and integrating the resulting equation over $\mathbb{R}^2$ lead to
\begin{align*}
& \int_{\mathbb{R}^2}\rho_0^R|\tilde{\mathbf{u}}_0^R|^2dx
+\mu\int_{\mathbb{R}^2}|\nabla\tilde{\mathbf{u}}_0^R|^2dx
+(\lambda+\mu)\int_{\mathbb{R}^2}|\divv\tilde{\mathbf{u}}_0^R|^2dx \\
& \leq \|\sqrt{\rho_0^R}\tilde{\mathbf{u}}_0^R\|_{L^2(B_R)}\|\mathbf{h}^R\|_{L^2(B_R)} +C\|P_0^R\|_{L^2(B_R)}\|\nabla\tilde{\mathbf{u}}^R_0\|_{L^2(B_R)}
+C\|\mathbf{b}_0^R\|_{L^4(B_R)}^2\|\nabla\tilde{\mathbf{u}}^R_0\|_{L^2(B_R)}
 \\
& \leq \frac12\int_{B_R}\rho_0^R |\tilde{\mathbf{u}}_0^R|^2dx
+\frac{\mu}{2}\|\nabla\tilde{\mathbf{u}}_0^R\|_{L^2(B_R)}^2+C\|\mathbf{h}^R\|_{L^2(B_R)}^2
+C\|P_0^R\|_{L^2(B_R)}^2+C\|\mathbf{b}_0^R\|_{L^4(B_R)}^4
 \\
& \leq\frac12\int_{B_R}\rho_0^R |\tilde{\mathbf{u}}_0^R|^2dx
+\frac{\mu}{2}\|\nabla\tilde{\mathbf{u}}_0^R\|_{L^2(B_R)}^2+C
\end{align*}
due to \eqref{4.2}--\eqref{4.4}, which implies
\begin{align}\label{4.7}
\int_{\mathbb{R}^2}\rho_0^R|\tilde{\mathbf{u}}_0^R|^2dx
+\int_{\mathbb{R}^2}|\nabla \tilde{\mathbf{u}}_0^R|^2dx \leq C
\end{align}
for some $C$ independent of $R.$
This together with \eqref{4.2} yields that there exist a subsequence $R_j\rightarrow \infty$ and a function $\tilde{\mathbf{u}}_0\in\{\tilde{\mathbf{u}}_0\in
H^1_{\loc}(\mathbb{R}^2)|\sqrt{\rho_0}\tilde{\mathbf{u}}_0\in L^2(\mathbb{R}^2), \nabla \tilde{\mathbf{u}}_0\in L^2(\mathbb{R}^2)\}$ such that
\begin{align}\label{4.8}
\begin{cases}
\sqrt{\rho^{R_j}_0}\tilde{\mathbf{u}}^{R_j}_0 \rightharpoonup \sqrt{\rho_0} \tilde{\mathbf{u}}_0 \mbox{ weakly in } L^2(\mathbb{R}^2) ,\\
\nabla \tilde{\mathbf{u}}_0^{R_j}\rightharpoonup \nabla \tilde{\mathbf{u}}_0 \mbox{ weakly in } L^2(\mathbb{R}^2).
\end{cases}
\end{align}
Next, we will show
\begin{align}\label{4.9}
\tilde{\mathbf{u}}_0=\mathbf{u}_0.
\end{align}
Indeed, subtracting \eqref{1.8} from \eqref{4.5} leads to
\begin{align}\label{4.10}
& -\mu\Delta\left(\tilde{\mathbf{u}}_0^{R_j}-\mathbf{u}_{0}\right)
-(\mu+\lambda)\nabla\divv\left(\tilde{\mathbf{u}}_0^{R_j}-\mathbf{u}_{0}\right)
+\nabla\left(P^{R_j}_0-P_0\right)
\notag \\ & =\mathbf{b}_0^{R_j}\cdot\nabla\mathbf{b}_0^{R_j}
-\mathbf{b}_0\cdot\nabla\mathbf{b}_0
-\frac12\nabla\left(|\mathbf{b}_0^{R_j}|^2-|\mathbf{b}_0|^2\right) \notag \\
& \quad +\left(\sqrt{\rho_{0}^{R_j}}\mathbf{g}*j_{1/R_j}-\sqrt{\rho_0}\mathbf{g}\right)
-\sqrt{\rho_{0}^{R_j}}\left(\sqrt{\rho_{0}^{R_j}}\tilde{\mathbf{u}}_0^{R_j}
-\mathbf{h}^{R_j}\right).
\end{align}
Multiplying \eqref{4.10} by a test function $\pmb{\phi}\in C_0^\infty(\mathbb{R}^2)$, it holds that
\begin{align}\label{4.11}
& \mu\int_{\mathbb{R}^2}\partial_i(\tilde{\mathbf{u}}_0^{R_j}-\mathbf{u}_0)
\cdot\partial_i\pmb{\phi}dx
+(\lambda+\mu)\int_{\mathbb{R}^2}\divv(\tilde{\mathbf{u}}_0^{R_j}-\mathbf{u}_0)
\cdot\divv\pmb{\phi}dx \notag \\
& \quad +\int_{\mathbb{R}^2} \sqrt{\rho^{R_j}_0}(\sqrt{\rho^{R_j}_0}\tilde{\mathbf{u}}^{R_j}_0-\mathbf{h}^{R_j})
\cdot\pmb{\phi}dx \notag \\
& = \int_{\mathbb{R}^2}\left(P^{R_j}_0-P_0\right)\divv\pmb{\phi}dx
+\int_{\mathbb{R}^2}\left(\sqrt{\rho_{0}^{R_j}}\mathbf{g}*j_{1/R_j}
-\sqrt{\rho_0}\mathbf{g}\right)\pmb{\phi}dx \notag \\
& \quad +\int_{\mathbb{R}^2}\left(\mathbf{b}_0^{R_j}-\mathbf{b}_0\right)
\cdot\nabla\mathbf{b}_0^{R_j}\cdot\pmb{\phi}dx
+\int_{\mathbb{R}^2}\mathbf{b}_0^{R_j}\cdot
\nabla\left(\mathbf{b}_0^{R_j}-\mathbf{b}_0\right)
\cdot\pmb{\phi}dx \notag \\
& \quad +\frac12\int_{\mathbb{R}^2}\left(|\mathbf{b}_0^{R_j}|^2-|\mathbf{b}_0|^2\right)
\divv\pmb{\phi}dx.
\end{align}
Let $R_j\rightarrow \infty$,
it follows from \eqref{4.2}, \eqref{4.3}, \eqref{4.4}, and \eqref{4.8} that
\begin{align}\label{4.12}
\int_{\mathbb{R}^2}\partial_i(\tilde{\mathbf{u}}_0
-\mathbf{u}_0)\cdot\partial_i\pmb{\phi}dx+\int_{\mathbb{R}^2}
\rho_0(\tilde{\mathbf{u}}_0-\mathbf{u}_0)\cdot\pmb{\phi}dx=0,
\end{align}
which implies \eqref{4.9}.
Furthermore, multiplying \eqref{4.5} by $\tilde{\mathbf{u}}_0^{R_j}$ and integrating the resulting equation over $\mathbb{R}^2$, by the same arguments as \eqref{4.12}, we have
\begin{align*}
\lim \limits_{R_j\rightarrow \infty}\int_{\mathbb{R}^2}
\left(|\nabla\tilde{\mathbf{u}}_0^{R_j}|^2
+\rho_0^{R_j}|\tilde{\mathbf{u}}_0^{R_j}|^2\right)dx
= \int_{\mathbb{R}^2}\left(|\nabla\mathbf{u}_0|^2+\rho_0|\mathbf{u}_0|^2\right)dx,
\end{align*}
which combined with \eqref{4.8} leads to
\begin{align*}
\lim\limits_{R_j\rightarrow \infty}\int_{\mathbb{R}^2}
|\nabla\tilde{\mathbf{u}}_0^{R_j} |^2dx
=\int_{\mathbb{R}^2}|\nabla\tilde{\mathbf{u}}_0|^2dx,\,\,
\lim\limits_{R_j\rightarrow\infty}\int_{\mathbb{R}^2}\rho_0^{R_j}   |\tilde{\mathbf{u}}_0^{R_j}|^2dx
=\int_{\mathbb{R}^2}\rho_0|\tilde{\mathbf{u}}_0|^2dx.
\end{align*}
This, along with  \eqref{4.9} and \eqref{4.8}, gives \eqref{4.6}.

Hence, by virtue of Lemma \ref{lem21}, the initial-boundary-value problem \eqref{2.2} with the initial data $(\rho_0^R,\mathbf{u}_0^R,P_0^R,\mathbf{b}_0^R)$ has a classical solution $(\rho^R,\mathbf{u}^R,P^R,\mathbf{b}^R)$ on $B_{R}\times [0,T_R]$. Moreover, Proposition \ref{pro} shows that there exists a $T_0$ independent of $R$ such that \eqref{3.3} holds for $(\rho^R,\mathbf{u}^R,P^R,\mathbf{b}^R)$.

For simplicity, in what follows, we denote
\begin{align*}
L^p=L^p(\mathbb{R}^2),\quad W^{k,p}=W^{k,p}(\mathbb{R}^2).
\end{align*}
Extending $(\rho^{R},\mathbf{u}^{R},P^{R},\mathbf{b}^{R})$ by zero on $\mathbb{R}^2\setminus B_{R}$ and denoting it by
\begin{align*}
\left(\tilde{\rho}^R\triangleq \varphi_R\rho^R,  \tilde{\mathbf{u}}^R, \tilde{P}^{R}, \tilde{\mathbf{b}}^R\right)
\end{align*}
with $\varphi_R$ as in \eqref{3.29}. We then infer from \eqref{3.3} that
\begin{align}\label{4.13}
& \sup\limits_{0\le t\le T_0}\left(\|\sqrt{\tilde{\rho}^R}\tilde{\mathbf{u}}^R\|_{L^2}
+\|\sqrt{\tilde{\rho}^R}\tilde{\dot{\mathbf{u}}}^R\|_{L^2}
+\|\nabla\tilde{\mathbf{u}}^R\|_{L^2}+\|\tilde{P}^R\|_{L^1\cap L^q}\right) \notag \\
& \leq \sup\limits_{0\le t\le T_0}\left(\|\sqrt{\rho^R}\mathbf{u}^R\|_{L^2}
+\|\sqrt{\rho^R}\dot{\mathbf{u}}^R\|_{L^2}\right)
+\sup\limits_{0\le t\le T_0}\left(\|\nabla\mathbf{u}^R\|_{L^2(B_R)}
+\|P^R\|_{L^1(B_R)\cap L^q(B_R)}\right) \notag \\
& \leq C
\end{align}
and
\begin{align}\label{4.14}
\sup\limits_{0\le t\le T_0}\|\tilde{\rho}^R\bar x^a\|_{L^1\cap L^\infty}\le C.
\end{align}
Similarly, it follows from \eqref{3.3}, \eqref{3.28}, and \eqref{3.38} that for $p\in[2,q]$,
\begin{align}\label{4.15}
\sup\limits_{0\leq t\leq T_0}
\left(\|\tilde{\mathbf{b}}^R\bar{x}^a\|_{L^p}
+\|\nabla\tilde{\mathbf{b}}^R\|_{L^p}\right)
+\int_0^{T_0}\|\nabla^2\tilde{\mathbf{u}}^R\|_{L^2\cap L^p}^2dt\leq C.
\end{align}
Next, for $p\in[2,q]$, we obtain from \eqref{3.3}, \eqref{3.28}, and \eqref{3.45} that
\begin{align}\label{4.16}
& \sup\limits_{0\le t\le T_0}
\left(\|\nabla(\tilde{\rho}^R\bar x^a)\|_{L^p}
+\|\nabla(\tilde{P}^R)\|_{L^p}\right) \notag \\
& \leq C\sup\limits_{0\le t\le T_0}\left( \|\nabla(\rho^R\bar x^a)\|_{L^p(B_R)}+R^{-1}\| \rho^R\bar x^a \|_{L^p(B_R)}\right) \notag \\
& \quad +C\sup\limits_{0\leq t\leq T_0}\left(\|\nabla(P^R)\|_{L^p(B_R)}+R^{-1}\|P^R\|_{L^p(B_R)}\right) \notag \\
& \leq C\sup\limits_{0\leq t\leq T_0}\left(\|\rho^R\bar x^a\|_{H^1(B_R)\cap W^{1,p}(B_R)}+\|P^R\|_{H^1(B_R)\cap W^{1,p}(B_R)}\right) \notag \\
& \leq C.
\end{align}

With the estimates \eqref{4.13}--\eqref{4.16} at hand, we find that the sequence
$(\tilde{\rho}^R,\tilde{\mathbf{u}}^R,\tilde{P}^R, \tilde{\mathbf{b}}^R)$ converges, up to the extraction of subsequences, to some limit $(\rho,\mathbf{u},P,\mathbf{b})$ in the
obvious weak sense, that is, as $R\rightarrow \infty,$ we have
\begin{align}\label{5.1}
& \tilde{\rho}^R\bar x\rightarrow \rho\bar{x}, \
\tilde{P}^R\rightarrow P, \
\mbox{in}\ C(\overline{B_N}\times [0,T_0]), \mbox{for any}N>0, \\
& \tilde{\rho}^R\bar x^a\rightharpoonup  \rho \bar x^a,\
\tilde{P}^R\rightharpoonup  P,\
\mbox{ weakly * in }L^\infty(0,T_0; H^1 \cap W^{1,q}), \\
& \tilde{\mathbf{b}}^R \bar{x}^a\rightharpoonup
\mathbf{b}\bar{x}^a,\ \nabla \tilde{\mathbf{b}}^R\rightharpoonup \nabla \mathbf{b},\ \mbox{weakly * in} \ L^\infty(0,T_0; L^2\cap L^q), \\
& \sqrt{\tilde{\rho}^R} \tilde{\mathbf{u}}^R\rightharpoonup \sqrt{\rho}\mathbf{u},
\,\, \nabla \tilde{\mathbf{u}}^R\rightharpoonup \nabla \mathbf{u},\,\,
\sqrt{\tilde{\rho}^R} \tilde{\dot{\mathbf{u}}}^R\rightharpoonup \sqrt{\rho}\dot{\mathbf{u}},
\mbox{ weakly * in }\ L^\infty(0,T_0; L^2), \\
& \nabla^2 \tilde{\mathbf{u}}^R\rightharpoonup \nabla^2\mathbf{u},\,\,
\nabla\tilde{P}^R\rightharpoonup  \nabla P, \mbox{ weakly  in}\ L^{2}(0,T_0;L^2\cap L^q), \end{align}
with
\begin{align}\label{5.2}
\rho \bar x^a\in L^\infty(0,T_0; L^1), \quad
\inf\limits_{0\leq t\leq T_0}\int_{B_{2N_0}}\rho(x,t)dx\geq \frac14.
\end{align}

Then letting $R\rightarrow \infty$, some standard arguments together with  \eqref{5.1}--\eqref{5.2} show that
$(\rho,\mathbf{u},P,\mathbf{b})$ is a strong solution of \eqref{mhd}-\eqref{n4} on $\mathbb{R}^2\times (0,T_0]$ satisfying \eqref{1.9} and \eqref{1.10}. The proof of Theorem 1.1 is completed.   \hfill $\Box$

\end{document}